\documentclass[a4paper, 12pt, twoside, notitlepage]{amsart}
\usepackage[utf8]{inputenc}
\usepackage[T1]{fontenc}
\usepackage{lmodern}
\usepackage[textwidth=14.5cm,textheight=22.5cm,centering]{geometry}
\usepackage{amsmath,amssymb,amsthm}
\usepackage{tikz}
\usetikzlibrary{arrows.meta}
\definecolor{figureblue}{RGB}{32,73,101}
\tikzset{surveyarrow/.style={-{Stealth[length=1.5mm]},line width=.5pt},
  surveycurve/.style={draw=figureblue,line width=.8pt}}
\usepackage[colorlinks=true,linkcolor=blue,citecolor=blue]{hyperref}

\newtheorem{theorem}{Theorem}

\newtheorem{corollary}[theorem]{Corollary}
\theoremstyle{definition}
\newtheorem{definition}[theorem]{Definition}
\theoremstyle{remark}

\newcommand{\R}{\mathbb{R}}
\newcommand{\WFa}{\operatorname{WF_{a}}}
\newcommand{\WF}{\operatorname{WF}}
\newcommand{\Cinf}{C^{\infty}}
\newcommand{\Comega}{C^{\omega}}

\begin{document}

\title{Double Fibration Transforms: Geometry and Microlocal Analysis}

\author{Shubham R. Jathar}
\address{Computational Engineering, School of Engineering Sciences,
Lappeenranta-Lahti University of Technology LUT, Lappeenranta, Finland}
\email{shubham.jathar@lut.fi}

\author[J. Railo]{Jesse Railo}
\address{Computational Engineering, School of Engineering Sciences,
Lappeenranta-Lahti University of Technology LUT, Lappeenranta, Finland}
\email{jesse.railo@lut.fi}


\begin{abstract}
Double fibration transforms integrate functions over families of submanifolds. We survey the recovery of singularities by both classical normal operator methods and direct analytic methods. We also discuss support theorems, injectivity and applications to geodesic ray transforms, seismic inversion and minimal surface data.
\end{abstract}
\maketitle
\section{Incidence geometry and examples}

The Radon transform asks whether a function is determined by integrals over hyperplanes \cite{Radon1917}. Double fibration transforms place this question in a geometric framework. An incidence relation between the object space and the data space specifies the integration manifolds. This viewpoint originates in Gelfand's program in integral geometry \cite{GGS1969,Helgason2011}. Guillemin and Sternberg developed the microlocal formulation \cite{Guillemin1985,GuilleminSternberg1979}.

We distinguish two recovery questions. Which singularities of the
unknown can be recovered stably from the data? When do vanishing
data force the unknown to vanish? Classical smooth microlocal
recovery uses the Bolker condition and ellipticity of the normal
operator. Mazzucchelli, Salo and Tzou instead prove analytic wave
front recovery directly from the transform, without using the
adjoint \cite[Section~1.4]{MST2023}. One-sided analytic continuation
then gives support results and, under additional geometric
hypotheses, injectivity.
\begin{definition}\label{def:dft}
Let $X$ and $G$ be smooth manifolds without boundary, of dimensions $n$ and $N$. Let $Z\subset G\times X$ be a closed embedded submanifold of codimension $k$, where $1\le k\le n-1$, whose projections
\[
\begin{tikzpicture}[baseline=(current bounding box.center),font=\small]
\node (Z) at (0,.65) {$Z$};
\node (G) at (-1.1,0) {$G$};
\node (X) at (1.1,0) {$X$};
\draw[surveyarrow] (Z) -- node[above left,inner sep=2pt] {$\pi_G$} (G);
\draw[surveyarrow] (Z) -- node[above right,inner sep=2pt] {$\pi_X$} (X);
\end{tikzpicture}
\]
are surjective submersions. For $z\in G$, we write $G_z=\{x\in X:(z,x)\in Z\}$ and $d=n-k$. Given a positive smooth weight $\kappa$ on $Z$ and positive fiber densities $d\mu_z$ varying smoothly with $z$, we define the associated \emph{double fibration transform} by
\[
Rf(z)=\int_{G_z}\kappa(z,x)f(x)\,d\mu_z(x),\qquad f\in\Cinf_c(X).
\]
\end{definition}

For each $z\in G$, the restriction $\pi_X|_{\pi_G^{-1}(z)}$ is a diffeomorphism onto the $d$-dimensional manifold $G_z$ \cite[proof of Lemma 2.1]{MST2023}. For $f\in\Cinf_c(X)$, the value $Rf(z)$ is the integral of the function $x\mapsto\kappa(z,x)f(x)$ over $G_z$ with respect to $d\mu_z$. After fixing smooth positive densities on $X$ and $G$, we define the formal adjoint by $$\langle Rf,h\rangle_G=\langle f,R^*h\rangle_X.$$ The densities defining $R$, together with the chosen densities on $X$ and $G$, induce the fiber densities for $R^*$ \cite[Section~1]{Quinto1980}. The submersion hypotheses and closedness of $Z$ give a continuous extension $R:\mathcal E'(X)\to\mathcal D'(G)$ \cite[Section~2]{MST2023}. The incidence relation determines the canonical relation, while the weighted fiber density $\kappa\,d\mu_z$ also enters the principal symbol \cite[Section~2]{Quinto1980}.

\begin{figure}[htbp]
\centering
\begin{tikzpicture}[
  font=\normalsize,
  x=1.2cm,y=1.2cm,
  surveyarrow/.append style={line width=.65pt,-{Stealth[length=2mm]}},
  surveycurve/.append style={line width=1pt}]
  \node at (2,1.5) {Data space $G$};
  \draw[surveyarrow] (-.2,0) -- (4.35,0) node[right] {$\theta$};
  \draw[surveyarrow] (0,-1) -- (0,1.18) node[above left] {$s$};
  \draw[surveycurve,domain=0:4,samples=90,smooth]
    plot (\x,{.6*cos(90*\x)+.5*sin(90*\x)});
  \draw[densely dashed,gray] (.666667,0) -- (.666667,.733013);
  \fill[figureblue] (.666667,.733013) circle (1.3pt);
  \node[above right,inner sep=2pt] at (.666667,.733013) {$z_0$};
  \node[below,inner sep=2pt] at (.666667,0) {$\theta_0$};
  \node[below,inner sep=2pt] at (2,0) {$\pi$};
  \node[below,inner sep=2pt] at (4,0) {$2\pi$};
  \node at (2,-1.18) {$s=x_0\cdot\omega(\theta)$};
  \begin{scope}[xshift=8.85cm]
    \node at (.5,1.5) {Object space $X$};
    \draw[surveyarrow,gray] (-1.05,0) -- (2,0) node[right,black] {$x_1$};
    \draw[surveyarrow,gray] (0,-.65) -- (0,1.15) node[left,black] {$x_2$};
    \draw[surveycurve] (-.85,1.337157) -- (1.85,-.221689);
    \node[below] at (1.7,-.221689) {$G_{z_0}$};
    \fill[figureblue] (.6,.5) circle (1.3pt);
    \node[below left,inner sep=2pt] at (.6,.5) {$x_0$};
    \draw[surveyarrow] (.6,.5) -- (.95,1.106218)
      node[right,inner sep=2pt] {$\omega(\theta_0)$};
    \node at (.5,-1.18) {$x\cdot\omega(\theta_0)=s_0$};
  \end{scope}
\end{tikzpicture}
\caption{Projected incidence fibers for planar lines. The curve on the left parameterizes lines through $x_0$. On the right, the point $z_0=(\theta_0,s_0)$ determines the line $G_{z_0}$ with the normal $\omega(\theta_0)$. Thus $(z_0,x_0)\in Z$.}
\label{fig:incidence}
\end{figure}
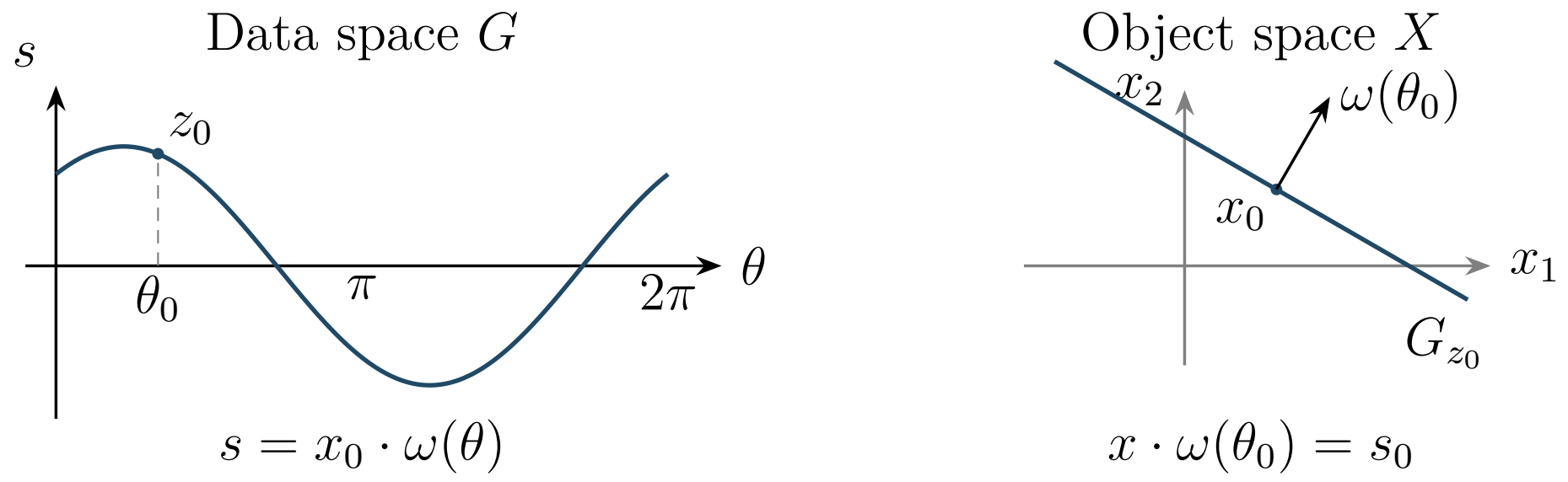

For hyperplanes in $\R^n$, we have $N=n$ and $d=n-1$. In the plane, oriented lines are parameterized by $(\theta,s)\in S^1\times\R$ as the solutions of the equation $x\cdot\omega(\theta)=s$, where $\omega(\theta)=(\cos\theta,\sin\theta)$. Figure~\ref{fig:incidence} illustrates the projected fibers. A fixed data point determines a line in $X$ and the collection of all lines through a fixed object point forms a curve in $G$. More generally, the family of affine $d$-planes has the dimension $N=(d+1)(n-d)$, which exceeds $n$ when $1\le d\le n-2$.

Geodesic ray transforms provide another class of examples. Let $(M,g)$ be a compact nontrapping Riemannian $n$-manifold with strictly convex boundary. Oriented maximal unit-speed geodesics are parameterized by strictly inward-pointing boundary unit vectors. Here $N=2n-2$ and $d=1$ \cite{PSU2023}. If all maximal geodesics are embedded, their incidence relation satisfies Definition~\ref{def:dft} with $X=M^{\mathrm{int}}$ \cite[Lemma~4.3 and Section~4.3]{MST2023}. Locally, fix an embedded maximal geodesic with initial vector $z_0$ and take $G$ to be a sufficiently small open neighborhood of $z_0$ in this boundary parameter space. The maximal geodesics with initial vectors in $G$ give a double fibration, with $X$ the open region covered by their interiors \cite[Proposition~4.7]{MST2023}. This embedded-incidence construction excludes
self-intersecting geodesics.

The geometric hypotheses matter for the injectivity of $R$. Outside the preceding nontrapping, strictly convex setting, the large spherical cap $M=\{x\in S^2:x_3\ge-\varepsilon\}$, with $0<\varepsilon<1$, contains trapped geodesics. A nonzero smooth function on $S^2$ that is odd under the antipodal map and supported in $\{|x_3|<\varepsilon\}$ integrates to zero over every boundary-to-boundary geodesic in $M$. The support restriction makes each such integral equal to the integral over the corresponding great circle, which vanishes by oddness \cite[Example~2.5.5]{PSU2023}. 

The light ray transform illustrates a different limitation. Timelike covectors annihilate no nonzero null tangent vector. Consequently, timelike singularities are microlocally invisible even with full data \cite[Section~3.1]{LOSU2020}. This obstruction to stable singularity recovery does not by itself imply nonuniqueness. On Minkowski space of spacetime dimension at least three, the unweighted light ray transform is injective on compactly supported distributions. The Fourier slice theorem determines the Fourier transform on the spacelike cone, and compact support permits analytic continuation to the remaining frequencies \cite[Section~3.1]{Stefanov2017}.

\section{The Bolker condition and its limits}

For $f\in\mathcal E'(X)$, the wave front set $\WF(f)$ records the positions and cotangent directions of singularities. The kernel of $R$ is conormal to $Z$. Thus $R$ is a Fourier integral operator (FIO) \cite[Theorem~2.2]{MST2023}. In H\"ormander's convention, the order and canonical relation are, respectively,
\[
\begin{aligned}
m&=\frac{k}{2}-\frac{N+n}{4}
  =-\frac d2+\frac{n-N}{4},\\
C&=(N^*Z\setminus0)'
  =\bigl\{(z,\zeta,x,\xi):
    (z,\zeta,x,-\xi)\in N^*Z\setminus0\bigr\}.
\end{aligned}
\]
Here $N^*Z$ is the conormal bundle of $Z$. The projections $\pi_L$ and $\pi_R$ from $C\subset(T^*G\setminus0)\times(T^*X\setminus0)$ select the data and object covectors, respectively. Since $C$ is Lagrangian, $d\pi_L|_c$ is injective if and only if $d\pi_R|_c$ is surjective, for each $c\in C$ \cite[Lemma~3.3]{MST2023}.

\begin{definition}\label{def:bolker}
The \emph{Bolker condition} holds at $c=(z,\zeta,x,\xi)\in C$ if $\pi_L^{-1}(z,\zeta)=\{c\}$ and $d\pi_L|_c$ is injective. It holds globally if $\pi_L$ is an injective immersion \cite[Definition~1.5(e) and Section~3]{MST2023}.
\end{definition}
Under the global Bolker condition, a data covector in $\pi_L(C)$ determines a unique object covector \cite[Definition~10]{Chihara2026}. Injectivity makes the composition $C^t\circ C$ diagonal, and the immersion condition gives the clean intersection required by the FIO composition calculus. To form the normal operator globally, we also assume that $\pi_X:Z\to X$ is proper \cite[proof of Theorem~13]{Chihara2026}. Without properness, support restrictions or localized transforms are needed.

\begin{theorem}[Guillemin
{\cite[Theorem~B.5.8]{StefanovUhlmannBook}}]
\label{thm:guillemin}
Let $R$ satisfy Definition~\ref{def:dft}. Assume that $\pi_X:Z\to X$ is proper and that $\pi_L$ is an injective immersion. Then $R^*R$ is an elliptic pseudodifferential operator of order $-d$ on $X$.
\end{theorem}
Retaining the geometric assumptions of Theorem~\ref{thm:guillemin}, we may allow the smooth weight $\kappa$ to vanish. With $C$ still defined by the incidence relation, we write
\[
\mathcal V_\kappa
=
\bigl\{(x,\xi)\in T^*X\setminus0:
\text{ there exists }(z,\zeta,x,\xi)\in C
\text{ with }\kappa(z,x)\ne0\bigr\}.
\]
The normal operator is elliptic on $\mathcal V_\kappa$, since a branch with nonzero weight contributes a positive term to the principal symbol \cite[Theorem~B.5.8]{StefanovUhlmannBook}. The ellipticity gives microlocal recovery of $f$ from $R^*Rf$ at these covectors.

Suppose that $K\subset X$ is compact and that $R^*R$ is elliptic at every nonzero covector based in $K$. Choose $\chi\in\Cinf_c(X)$ equal to one near $K$. A localized elliptic parametrix and the Sobolev mapping properties of pseudodifferential operators \cite[Sections~A.5.12--A.5.13]{StefanovUhlmannBook} give, for every $s\in\R$, the estimate
\[
\|f\|_{H^s}
\le C\bigl(
\|\chi R^*Rf\|_{H^{s+d}}
+\|f\|_{H^{s-1}}
\bigr),
\qquad
f\in H^s(X),\quad \operatorname{supp}f\subset K.
\] Here $C$ depends on $R$, $K$, $\chi$, and $s$.
If $R$ is injective on distributions supported in $K$, the term $\|f\|_{H^{s-1}}$ can be removed by a compactness argument \cite[Lemma~IV.2.1]{StefanovUhlmannBook}.

Conjugate points can obstruct the stable recovery of singularities from geodesic X-ray data. On Riemannian surfaces, conjugate points can destroy injectivity of $\pi_L$ even though $\pi_L$ remains a local diffeomorphism \cite[Theorem~3.2]{MonardStefanovUhlmann2015}. Under the fold-caustic hypotheses of Stefanov and Uhlmann, the localized normal operator is the sum of a pseudodifferential operator of order $-1$ and an FIO of order $-n/2$, where $n=\dim M$ \cite[Theorem~2.1]{StefanovUhlmann2012}. These orders coincide when $n=2$, and singularities at conjugate points can cancel for the unweighted geodesic X-ray transform, obstructing stable recovery \cite[Corollary~4.1 and Theorem~4.4(b)]{MonardStefanovUhlmann2015}. In the fold-caustic setting with $n\ge3$, a local canonical graph condition on the FIO term gives microlocal invertibility \cite[Corollary~9.1]{StefanovUhlmann2012}.

More general decompositions require hypotheses on the conjugate points. For geodesic X-ray transforms with regular conjugate pairs, Holman and Uhlmann obtain FIO terms indexed by conjugate-pair multiplicity under their geometric hypotheses \cite[Theorem~4]{HolmanUhlmann2018}. Chihara proves an analogous normal operator decomposition for double fibration transforms with regular $Z$-conjugate triplets under additional hypotheses \cite[Theorem~14]{Chihara2026}.

Even when the global Bolker condition excludes this geometric obstruction, the transform may annihilate nonzero smooth functions. Boman's example separates singularity recovery from injectivity.

\begin{theorem}[Boman {\cite{Boman1993}}]\label{thm:boman}
There exist a strictly positive smooth planar line weight $\kappa$ and a nonzero $f\in\Cinf_c(\R^2)$ such that $R_\kappa f=0$ on every oriented line.
\end{theorem}

The planar point-line incidence relation satisfies the global Bolker condition, and the projection to the object space is proper. Theorem~\ref{thm:guillemin} therefore makes $R_\kappa^*R_\kappa$ elliptic of order $-1$. Ellipticity recovers singularities but leaves open the possibility of nonzero smooth functions in the kernel.

\section{Analytic support and injectivity}

Analytic microlocal analysis provides a route from singularity recovery to uniqueness. On a real-analytic manifold $X$, the analytic wave front set $\WFa(f)$ detects failure of analytic regularity in cotangent directions. If a distribution $f$ vanishes on one side of a $C^2$ hypersurface $\Sigma$ near a point $x_0\in\operatorname{supp}f$, both conormal directions to $\Sigma$ at $x_0$ belong to $\WFa(f)$ \cite[Theorem~6.1]{MST2023}. Excluding either direction therefore forces local vanishing. This mechanism underlies analytic support theorems for hypersurface transforms \cite{BomanQuinto1987}, geodesic rays \cite{Krishnan2009} and light rays \cite{Stefanov2017}.

\begin{theorem}[Mazzucchelli--Salo--Tzou {\cite[Theorem~1.2]{MST2023}}]\label{thm:mst}
Let $R$ satisfy Definition~\ref{def:dft}, with $N\ge n$. Assume that the manifolds, incidence relation, fiber densities and nonvanishing weight are real-analytic. If the Bolker condition holds at $(z,\zeta,x,\xi)\in C$, then, for every $f\in\mathcal E'(X)$, we have
\[
(z,\zeta)\notin\WFa(Rf)\ \Longrightarrow\ (x,\xi)\notin\WFa(f).
\]
\end{theorem}

The proof works directly with $R$. Analytic wave packets transfer analytic regularity from the data covector to the unique object covector selected by the Bolker condition \cite[Section~1.4]{MST2023}. Locally, the analytic fiber densities differ from the orientation forms in the source by a nonvanishing analytic factor, which can be absorbed into the weight. Combining the theorem with one-sided continuation gives the following support result.

\begin{figure}[htbp]
\centering
\begin{tikzpicture}[
  font=\normalsize,
  x=1.2cm,y=1.2cm,
  surveyarrow/.append style={line width=.65pt,-{Stealth[length=2mm]}},
  surveycurve/.append style={line width=1pt}]
  \fill[black!7] plot[domain=-1.8:1.8,samples=50] (\x,{.3*\x*\x})
    -- (1.8,1.35) -- (-1.8,1.35) -- cycle;
  \draw[surveycurve,domain=-1.8:1.8,samples=50,smooth]
    plot (\x,{.3*\x*\x});
  \node[right,inner sep=2pt] at (1.8,.972) {$\Sigma$};
  \draw[densely dashed,line width=.6pt] (-1.95,0) -- (2.05,0)
    node[right,inner sep=2pt] {$G_{z_0}$};
  \fill (0,0) circle (1.3pt);
  \node[below left,inner sep=2pt] at (0,0) {$x_0$};
  \draw[surveyarrow] (0,0) -- (0,.65) node[right,inner sep=2pt] {$\xi_0$};
  \node at (0,1.1) {unknown};
  \node at (-.6,-.5) {$f=0$};
  \begin{scope}[xshift=7.65cm]
    \node at (0,1.4) {$Rf=0$ near $z_0$};
    \draw[surveyarrow] (0,1.07) -- (0,.68);
    \node at (0,.35) {$(x_0,\xi_0)\notin\WFa(f)$};
    \draw[surveyarrow] (0,.02) -- (0,-.37);
    \node at (0,-.70) {$f=0$ near $x_0$};
  \end{scope}
\end{tikzpicture}
\caption{Local support under analytic and pointwise Bolker hypotheses. Vanishing data exclude an analytic conormal singularity. One-sided unique continuation then forces local vanishing. A tangent integration curve is shown schematically.}
\label{fig:continuation}
\end{figure}
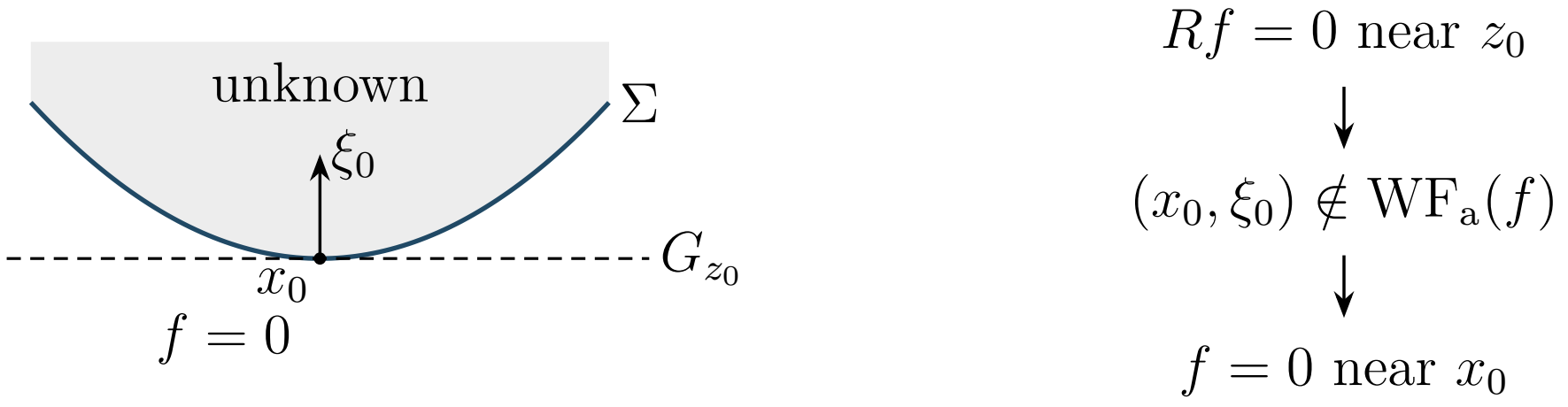

\begin{corollary}[Local support {\cite[Theorem~1.3]{MST2023}}]\label{cor:local-support}
Under the analytic assumptions of Theorem~\ref{thm:mst}, let $f\in\mathcal E'(X)$ vanish on one side of a $C^2$ hypersurface $\Sigma$ near $x_0$. Suppose that $0\ne\xi_0\in N^*_{x_0}\Sigma$ and that the Bolker condition holds at $(z_0,\zeta_0,x_0,\xi_0)\in C$. If $Rf=0$ near $z_0$, then $f=0$ near $x_0$.
\end{corollary}
Figure~\ref{fig:continuation} illustrates the two implications. Vanishing data exclude the analytic conormal singularity that would occur at a boundary point of the support. Therefore, the support cannot meet the hypersurface at that point.

Chihara, Jathar and Railo extend Theorem~\ref{thm:mst} to
vector-valued distributions and invertible analytic matrix
weights \cite[Theorem~1.2]{CJR2026}. Their local result for the geodesic ray transform assumes
a compact real-analytic Riemannian manifold with boundary
and an invertible real-analytic matrix weight. Short geodesics near a strictly convex boundary point give
local injectivity in every dimension $n\ge2$
\cite[Theorem~1.5]{CJR2026}. Compact real-analytic Riemannian surfaces with strictly convex boundary and nontrapping geodesic flow admit a strictly convex foliation. Repeated application of the local injectivity theorem along this foliation yields global injectivity
\cite[Corollary~1.6]{CJR2026}. The scalar case is established in \cite[Theorem~1.4]{MST2023}.

\begin{theorem}[Analytic surfaces {\cite[Corollary~1.6]{CJR2026}}]\label{thm:analytic-surface}
Let $(M,g)$ be a compact real-analytic Riemannian surface with strictly convex boundary. Assume that the geodesic flow is nontrapping. For an integer $r\ge1$, let $W\in\Comega(SM;\mathrm{GL}(r,\mathbb C))$, where $SM$ is the unit sphere bundle. For each oriented maximal unit-speed geodesic $\gamma$ of length $\ell_\gamma$, define the matrix weighted geodesic ray transform as
\[
I_Wf(\gamma)=\int_0^{\ell_\gamma}
W(\gamma(t),\dot\gamma(t))f(\gamma(t))\,dt.
\]
Then $I_W$ is injective on $\Cinf(M;\mathbb C^r)$.
\end{theorem}
Conjugate points are allowed, and the unknown need not be analytic. The proof uses short geodesic segments near the leaves of the foliation. No global Bolker condition is required.

On smooth simple surfaces, Paternain and Salo determine general complex connections and Higgs fields from full nonabelian data up to a gauge transformation equal to the identity on the boundary \cite[Theorem~1.1]{PS2024}. The associated linear theorem concerns attenuations defined by connections and Higgs fields \cite[Theorem~1.2]{PS2024}. These results do not assert injectivity for arbitrary smooth matrix weights.

In dimensions $n\ge3$, Uhlmann and Vasy use the scattering calculus to obtain stable local recovery of functions near a strictly convex boundary point and global injectivity under suitable foliations by strictly convex hypersurfaces \cite{UhlmannVasy2016}. Paternain, Salo, Uhlmann and Zhou extend local recovery to arbitrary smooth invertible matrix weights \cite[Theorem~1.3]{PSUZ2019}. Applying the layer-stripping argument of \cite[Section~6]{PSUZ2019} to this local theorem gives global injectivity on smooth vector-valued functions when the compact Riemannian manifold has strictly convex boundary and admits a smooth strictly convex function. The analytic local uniqueness theorem applies in every
dimension $n\ge2$ \cite[Theorem~1.5]{CJR2026},
whereas these smooth local recovery theorems require $n\ge3$.

\section{Recent applications and questions for the future}

The same analytic recovery mechanism applies to linearizations that are more general FIOs. Busch proves analytic wave front recovery for a class of elliptic classical analytic FIOs with frequency-dependent amplitudes, under the Bolker condition and additional phase and canonical relation hypotheses \cite[Theorem~3.1 and Corollary~3.2]{Busch2025}. An analytic hyperbolic parametrix represents a linearized acoustic seismic operator as an FIO \cite[Theorems~4.1 and~5.4]{Busch2025}. For an analytic background wave speed and the prescribed geometric and acquisition assumptions, analytic microlocal continuation and layer stripping give injectivity on continuous perturbations supported in a compact subsurface region \cite[Theorem~5.1]{Busch2025}.

For a minimal-hypersurface area data, Busch, Liimatainen, Salo and Tzou obtain a scalar minimal surface transform by conformal linearization \cite[Proposition~2.3]{BLST2025}. Suitable finite-dimensional families satisfy the Bolker condition \cite[Corollary~3.3]{BLST2025}.
Under the hypothesis of \cite[Theorem~4.2]{BLST2025}, the authors prove local uniqueness and conditional H\"older
stability for small smooth conformal perturbations of an
analytic background in ambient dimension at least three. 

The passage from microlocal recovery to uniqueness raises three questions. \textbf{Q1:} Which regularity classes between $\Comega$ and $\Cinf$ for the geometry, weight and fiber densities admit support theorems under the Bolker condition? \textbf{Q2:} Is the unweighted geodesic ray transform injective on smooth functions on every compact nontrapping Riemannian surface with strictly convex boundary, extending the analytic result of \cite[Corollary~1.6]{CJR2026}? \textbf{Q3:} Which choices of weights or additional measurements permit stable separation when several conjugate branches contribute to the same data covector? Boman's example shows that the Bolker condition and smooth positive weights do not guarantee injectivity, while conjugate-point cancellation can obstruct stable recovery.

\section*{Acknowledgements}
The work of S.R.J. and J.R. was supported by the Research Council of Finland through the Flagship of Advanced Mathematics for Sensing, Imaging and Modelling (decision number 359183) and the Jenny and Antti Wihuri Foundation. In addition, J.R. was supported by the Emil Aaltonen Foundation.

\end{document}